\documentclass[a4paper,twoside,reqno]{amsart}
\usepackage{xypic,amscd,amssymb,latexsym,srcltx}
\usepackage{dsfont}
\usepackage{amsthm}
\def\@setcopyright{}
\makeatother
\usepackage{mathtools}

\renewcommand{\a}{\underline{\mathrm a}}

\newtheorem*{acknowledgments*}{ACKNOWLEDGMENTS}

\newtheorem*{def*}{Definition}
\newtheorem*{thm*}{Wedderburn-Artin Theorem}
\newtheorem*{remark*}{Remark}
\newtheorem*{lemma*}{Lemma 1}
\newtheorem*{lemma**}{Lemma 2}
\newtheorem*{example*}{[Example]}

\begin{document}
\begin{center}
{\Large \bf A Constructive Proof of the Wedderburn-Artin Theorem
}\\
\bigskip
\end{center}
 \begin{center}
Sheng Gao
\footnote{Foundation item: The work is supported by National Natural Science Foundation of China(No.11701134).}
\footnote{E-mail: highlyrise@126.com}

\end{center}
\begin{center}
 \footnotesize{\emph{School of Mathematics, Hefei University of Technology, Hefei 230000, China }}

\end{center}

\title{}
\def\abstractname{\textbf{Abstract}}

\begin{abstract}\addcontentsline{toc}{section}{\bf{English Abstract}}
  In this paper, we use the idempotent decomposition to  give
  an explicit isomorphism from an arbitrary semisimple Artinian ring to an external
direct sum of finitely many full matrix rings over division rings.

\end{abstract}

\maketitle
The Wedderburn-Artin theorem on the structure of semisimple Artinian rings is one of the most fundamental results in the representation theory of groups and algebras. One can refer to \cite{Artin1}, \cite{Artin2}, \cite{Artin3} and \cite{Wedderburn} for its original version  and  early generalizations. Since the establishment of this theorem, new proofs  have been continuously proposed, most of which involved the properties of idempotent decompositions  and completely reducible modules.

In 1954, T.Szele \cite{T.Szele} gave a proof of the Wedderburn-Artin theorem on the basis of the general density theorem of  Chevalley and Jacobson. In 1965, D. W. Henderson \cite{D.W.Henderson} provided a simpler module-theoretic  proof by means of a key lemma: If $A$ is a ring with identity $1_{A}$ and  with an idempotent $e$ such that $AeA=A$, then $A$ is isomorphic to the endomorphism ring of $(Ae)_{eAe}$.  In 1993, W.K.Nicholson \cite{W.K.Nicholson} used a long-neglected lemma of Brauer on minimal left ideals to prove a semiprime version of the Wedderburn-Artin theorem.  In recent years,  M.Bre\v{s}ar (\cite{M.Bresar1} and \cite{M.Bresar2}) carried out a systematic  study on the structure of semiprime Artinian rings, and obtained some fascinating results similar to \cite{W.K.Nicholson}. In 2017, Tsiu-Kwen Lee \cite{Tsiu-Kwen Lee} presented another elegant proof of the classical version of this theorem, which  works for rings that are not necessarily unital.

Unfortunately,  in the existing literature we cannot find  the specific isomorphism from an arbitrary semisimple Artinian ring to an external direct sum of finitely many full matrix rings over  division rings. In  this paper, we will construct this isomorphism using the orthogonal decomposition of idempotents as a tool. It is worth noting that the method used in this paper is similar to \cite{M.Bresar1} and \cite{M.Bresar2} (e.g.  the application of matrix units), but the proof we give here is more constructive.

Throughout this paper, we let  $A$ denote a ring with identity and let $M_{n}(A)$ denote the full matrix ring of degree $n$ over $A$. A $n\times n$-matrix whose $(i,j)$-entry is $a_{ij}$ is abbreviated  to  $\big(a_{ij}\big)_{n\times n}$. The intersection of all maximal right ideals of $A$ is called the Jacobson radical of $A$ and is denoted by $J(A)$. $A$ is said to be semisimple if $J(A)=0$. A right $A$-module $V$ is denoted by $V_{A}$. $V_{A}$ is said to be Artinian if the lattice of all submodules of $V$ satisfies the descending chain condition, i.e., for every infinite chain of $A$-submodules of $V$
\begin{equation*}
V_{1}\supseteq V_{2} \supseteq V_{3}\supseteq \cdots,
\end{equation*}
there exists a number $m$ such that $V_{m}=V_{m+1}=V_{m+2}=\cdots$. $A$ is said to be a right Artinian ring if the right regular module $A_{A}$ is an Artinian module (see \cite[Chapter 3]{Anderson}). It is easy to show that $J(A)$ is a two-sided ideal of $A$ and coincides with the intersection of all maximal left ideals of $A$ (see \cite[Chapter 1]{Hirosi Nagao}).
\begin{def*}
 If $e$ and $f$ are two idempotents of $A$ such that $eA$ and  $fA$ are isomorphic as right $A$-module, then  $e$ and $f$ are  said to be equivalent and we write $e\simeq f$.
\end{def*}
Concerning the idempotents, all the facts needed for our discussion are collected in the following lemma.
\begin{lemma*}\cite[Theorem 4.1-4.4, pp.16-18]{Hirosi Nagao}
Let $e$ and $f$ be idempotents of $A$.
\item[\rm(1)] If $(eA)_{A}$ is an internal direct sum of $A$-submodules $I_{i}$$(1\leq i \leq n)$:
\[eA=I_{1}\oplus I_{2}\oplus\cdots\oplus I_{n},\]
then there exists an orthogonal idempotent decomposition \[e=e_{1}+e_{2}+\cdots+e_{n}\] such that
$I_{i}=e_{i}A\ (1\leq i\leq n)$. $e$ is  primitive if and only if $(eA)_{A}$ is indecomposable. In particular, if $(eA)_{A}$ is irreducible, then $e$ is  primitive.
\item[\rm(2)]$e\simeq f$ if and only if there exist $a\in fAe, b\in eAf$
such that $ab=f,ba=e$.
\item[\rm(3)]$\mathrm{Hom}_{A}(eA,fA)\cong fAe$ as additive groups and
$\mathrm{End}_{A}(eA)\cong eAe$ as rings.
\end{lemma*}

\begin{lemma**}\cite[Lemma 1]{M.Bresar1}
If $A$ contains a family elements $\varepsilon_{\mu\nu}(1\leq \mu,\nu\leq n)$ such that
$\varepsilon_{11}+\varepsilon_{22}+\cdots+\varepsilon_{nn}=1_{A}$ and $\varepsilon_{\mu\nu}\varepsilon_{\xi\eta}=\delta_{\nu\xi}\varepsilon_{\mu\eta}$
for all $\mu,\nu,\xi,\eta\in\{1,2,\cdots,n\}$$($here $\delta_{\nu\xi}$ is the Kronecker symbol$)$, then we have the following ring isomorphism
 \begin{eqnarray*}
 \varphi:A&\longrightarrow & M_{n}(\varepsilon_{11}A\varepsilon_{11}).\\
     a&\longmapsto&(\varepsilon_{1\mu}a\varepsilon_{\nu1})_{n\times n}=\left(\begin{matrix}\varepsilon_{11}a\varepsilon_{11}& \varepsilon_{11}a\varepsilon_{21}&\cdots&\varepsilon_{11}a\varepsilon_{n1}\\\varepsilon_{12}a\varepsilon_{11}& \varepsilon_{12}a\varepsilon_{21}&\cdots&\varepsilon_{12}a\varepsilon_{n1}\\\cdots&\cdots&\cdots&\cdots\\\varepsilon_{1n}a\varepsilon_{11}& \varepsilon_{1n}a\varepsilon_{21}&\cdots&\varepsilon_{1n}a\varepsilon_{n1}\end{matrix}\right)
\end{eqnarray*}
The elements $\varepsilon_{\mu\nu}$'s are called matrix units in $A$.
\end{lemma**}
Now, we will  prove the main result of this article.
\begin{thm*}
If $A$ is a semisimple right Artinian ring, then $A$ is isomorphic to an external
direct sum of finite full matrix rings over division rings.
\end{thm*}
\begin{proof}
 It follows from \cite[Theorem 8.1, pp.31]{Hirosi Nagao} that the right regular module $A_{A}$ is completely reducible, so we can express $A_{A}$ as an internal direct sum of a family of irreducible $A$-submodules $\{I_{\lambda}\}_{\lambda\in \Lambda}$:
$A=\bigoplus_{\lambda\in\Lambda}I_{\lambda}$. Then there exist finitely many elements $\lambda_{1},\lambda_{2},\cdots,\lambda_{n}\in\Lambda$
 such that $1_{A}\in \bigoplus_{i=1}^{n}I_{\lambda_{i}}$. This yields that $A=\bigoplus_{i=1}^{n}I_{\lambda_{i}}$ and $\Lambda=\{\lambda_{1},\lambda_{2},\ldots,\lambda_{n}\}$. Thus it follows from Lemma 1(1) that $1_{A}$ has a orthogonal primitive idempotent decomposition: \[1_{A}=\sum_{i=1}^{k}\sum_{\mu=1}^{n_{i}}e_{i\mu},\]
  where $\{I_{\lambda_{1}},I_{\lambda_{2}},\cdots,I_{\lambda_{n}}\}=\{e_{i\mu}A|1\leq i \leq k, 1 \leq \mu \leq n_{i}\}$ and  $e_{i\mu}$'s
are pairwise orthogonal primitive idempotents of $A$ such that $e_{i\mu}\simeq e_{j\nu}$ iff $\ i=j$.

For each $i\in \{1,2,\cdots,k\}$  and  $\mu\in \{1,2,\cdots,n_{i}\}$, we can choose $a_{\mu}^{(i)}\in e_{i1}Ae_{i\mu}$
and  $b_{\mu}^{(i)}\in e_{i\mu}Ae_{i1}$ such that $e_{i1}=a_{\mu}^{(i)}b_{\mu}^{(i)}$ and $e_{i\mu}=b_{\mu}^{(i)}a_{\mu}^{(i)}$ by Lemma 1(2). In particular, we set $a_{1}^{(i)}=b_{1}^{(i)}=e_{i1}$. According to Schur's lemma and Lemma 1(3), $e_{i1}Ae_{i1}\cong \mathrm{End}_{A}(e_{i1}A)$ is a division ring
and  $e_{i\mu}A e_{j\nu}\cong \mathrm{Hom}_{A}(e_{j\nu}A,e_{i\mu}A )=0$ whenever $i\neq j$.
Furthermore, we can show that
\begin{equation*}
a_{\nu}^{(i)}b_{\xi}^{(i)}=\delta_{\nu\xi}e_{i1} \tag{$\clubsuit$}
\end{equation*}
for each $i\in \{1,2,\cdots,k\}$  and
 $\nu,\xi\in \{1,2,\cdots,n_{i}\}$. This is evidently true when $\nu = \xi$. If $\nu \neq \xi$,
then it follows from $e_{i\nu}e_{i\xi}=0$ that
\begin{equation*}
 a_{\nu}^{(i)}b_{\xi}^{(i)}\in (e_{i1}Ae_{i\nu})(e_{i\xi}Ae_{i1})= e_{i1}A(e_{i\nu}e_{i\xi})Ae_{i1}=\{0\}.
\end{equation*}

For each $i\in \{1,2,\cdots,k\}$, set
$c_{i}=\sum_{\mu=1}^{n_{i}}e_{i\mu}$. Then we have $1_{A}=\sum_{i=1}^{k}c_{i}$ and
\begin{eqnarray*}
c_{i}Ac_{j}=\left(\sum_{\mu=1}^{n_{i}}e_{i\mu}\right)A\left(\sum_{\nu=1}^{n_{j}}e_{j\nu}\right)=\sum_{\mu=1}^{n_{i}}\sum_{\nu=1}^{n_{j}}e_{i\mu}A e_{j\nu}=0
\end{eqnarray*}
whenever $i\neq j$, which yields that
$c_{i}^{2}=\sum_{j=1}^{k}c_{i}c_{j}=c_{i}\sum_{j=1}^{k}c_{j}=c_{i}1_{A}=c_{i}$ and
\begin{equation*}
 c_{i}a=c_{i}a\sum_{j=1}^{k}c_{j}=\sum_{j=1}^{k}c_{i}ac_{j}=c_{i}ac_{i}=\sum_{j=1}^{k}c_{j}ac_{i}=\left(\sum_{j=1}^{k}c_{j}\right)ac_{i}=ac_{i}
\end{equation*}
for each $i\in\{1,2,\cdots,k\}$ and $a\in A$. Hence $c_{i}$'s are all central idempotents of $A$ and they are pairwise orthogonal. From this fact  we conclude that
$A=\bigoplus_{i=1}^{k}Ac_{i}$, where each $Ac_{i}$ is a two-sided ideal of $A$.

Set $A_{i}=Ac_{i}=c_{i}A \ (1\leq i\leq k)$. For each $i\in\{1,2,\cdots,k\}$, $A_{i}$ is obviously a ring with identity $c_{i}$. Using Lemma 2, we  can easily construct an isomorphism  from $A_{i}$  to the $n_{i}\times n_{i}$ full matrix ring over $e_{i1}Ae_{i1}$. First, we have
\begin{equation*}
  e_{i\mu}=\sum_{\nu=1}^{n_{i}}e_{i\mu}e_{i\nu}=e_{i\mu}\sum_{\nu=1}^{n_{i}}e_{i\nu}=e_{i\mu}c_{i}=c_{i}e_{i\mu}\in A_{i}
\end{equation*}
 for each $\mu\in \{1,2,\cdots,n_{i}\}$, so
 $ b_{\mu}^{(i)}, a_{\nu}^{(i)}\in A_{i}$ for each $\mu,\nu $ $(1\leq \mu,\nu \leq n_{i})$. Then
an easy  computation shows that
\begin{equation*}
\sum_{\mu=1}^{n_{i}}b_{\mu}^{(i)}a_{\mu}^{(i)}=\sum_{\mu=1}^{n_{i}}e_{i\mu}=c_{i}=1_{A_{i}}.
\end{equation*}
Moreover, by $(\clubsuit)$  we see that
\begin{equation*}
(b_{\mu}^{(i)}a_{\nu}^{(i)})(b_{\xi}^{(i)}a_{\eta}^{(i)})=b_{\mu}^{(i)}(a_{\nu}^{(i)}b_{\xi}^{(i)})a_{\eta}^{(i)}=b_{\mu}^{(i)}(\delta_{\nu\xi}e_{i1})a_{\eta}^{(i)}=\delta_{\nu\xi}b_{\mu}^{(i)}a_{\eta}^{(i)}
\end{equation*}
for each $\mu,\nu,\eta,\xi\in\{1,2,\cdots,n_{i}\}$.
Therefore, by Lemma 2, we assert that these elements $b_{\mu}^{(i)}a_{\nu}^{(i)}(1\leq\mu,\nu\leq n_{i})$ are matrix units in $A_{i}$ , which yields the following isomorphism of rings:
\begin{eqnarray*}
\Upsilon_{i}:A_{i}&\longrightarrow &M_{n_{i}}(b_{1}^{(i)}a_{1}^{(i)}A_{i}b_{1}^{(i)}a_{1}^{(i)})=M_{n_{i}}(e_{i1}Ae_{i1}),\\
ac_{i}&\longmapsto& \big(b_{1}^{(i)}a_{\mu}^{(i)}ac_{i}b_{\nu}^{(i)}a_{1}^{(i)}\big)_{n_{i}\times n_{i}}=\big(a_{\mu}^{(i)}ab_{\nu}^{(i)}\big)_{n_{i}\times n_{i}}\\
(\forall a \in A)& &\ \
\end{eqnarray*}
where  $\big(a_{\mu}^{(i)}ab_{\nu}^{(i)}\big)_{n_{i}\times n_{i}}$ denotes a $n_{i}\times n_{i}$-matrix whose $(\mu,\nu)$-entry is $a_{\mu}^{(i)}ab_{\nu}^{(i)}\ (1\leq\mu,\nu\leq n_{i})$.

From the above, we obtain an isomorphism of rings
\begin{eqnarray*}
\Upsilon:A&\longrightarrow &  M_{n_{1}}(e_{11}Ae_{11})\bigoplus M_{n_{2}}(e_{21}Ae_{21})\bigoplus \cdots\bigoplus M_{n_{k}}(e_{k1}Ae_{k1})\ ,\\
\end{eqnarray*}
which is defined by
\begin{eqnarray*}
&&\Upsilon(a)\\
&=&\big(\Upsilon_{1}(ac_{1}),\Upsilon_{2}(ac_{2}),\cdots,\Upsilon_{k}(ac_{k})\big)\\
&=&\left(\big(a_{\mu}^{(1)}ab_{\nu}^{(1)}\big)_{n_{1}\times n_{1}},\big(a_{\mu}^{(2)}ab_{\nu}^{(2)}\big)_{n_{2}\times n_{2}},\cdots,\big(a_{\mu}^{(k)}ab_{\nu}^{(k)}\big)_{n_{k}\times n_{k}}\right)\\
&=& \left(\left(\footnotesize\begin{matrix}a_{1}^{(1)}a b_{1}^{(1)}&\cdots&a_{1}^{(1)}a b_{n_{1}}^{(1)}\\\cdots&\cdots&\cdots\\a_{n_{1}}^{(1)}a b_{1}^{(1)} &\cdots &a_{n_{1}}^{(1)}a b_{n_{1}}^{(1)}\end{matrix}\right),\footnotesize\left(\begin{matrix}a_{1}^{(2)}a b_{1}^{(2)} &\cdots &a_{1}^{(2)}a b_{n_{2}}^{(2)}\\\cdots&\cdots&\cdots\\a_{n_{2}}^{(2)}a b_{1}^{(2)} &\cdots &a_{n_{2}}^{(2)}a b_{n_{2}}^{(2)}\end{matrix}\right),\cdots,\footnotesize\left(\begin{matrix}a_{1}^{(k)}a b_{1}^{(k)} &\cdots &a_{1}^{(k)}a b_{n_{k}}^{(k)}\\\cdots&\cdots&\cdots\\a_{n_{k}}^{(k)}a b_{1}^{(k)} &\cdots &a_{n_{k}}^{(k)}a b_{n_{k}}^{(k)}\end{matrix}\right)\right)
\end{eqnarray*}
for any $a\in A$. This completes the proof.
\end{proof}

\textbf{ACKNOWLEDGMENTS}\hfil\break
Thanks are given to Prof.Nicholson  for his helpful comments and suggestions.

\end{document}